# Multi-Agent Interaction in Social Trading Network

Oleg Malafeyev[1a], Nadezhda Redinskikh[2a], Nikolay Rumyantsev[3a]

[a] Saint-Petersburg State University,
Faculty of Applied Mathematics and Control Processes
Department of Modelling in Social and
Economical Systems
Saint-Petersburg, Russia

**Abstract**

The online retailers network models are considered. In some nodes of the network consumers are located. Each consumer wishes to purchase a particular product at minimal cost due to the price of goods and transport corruption costs. Also, in some nodes of network online-retailers wish to allocate the points of goods delivery. The points of goods delivery must be placed in accordance with a certain principle of optimality. In this paper we propose an algorithm for finding the optimal placement of goods delivery points in accordance with the compromise solution as the principle of optimality.

## 1. INTRODUCTION

E-Commerce has become an integral part of the economy of any state. If several years ago people were apprehensive about buying goods via the Internet, it is now hard to find anyone who has never taken advantage of an online store. The number of operations in this area increases incredibly fast, and E-commerce is having more and more impact on the global economy.

The network game model is considered in [15]. The model of corruption networks based on electric-network analogues is formalized and studied in this paper. Algorithms for finding the Cournot-Nash equilibrium [1] and the compromise point are proposed in this model. Other game-theoretic models, that can be applied to the online retail are considered in [2-25]. The online-retailers who sell goods over the Internet shops are considered in this paper. An online retailer needs to build the points of goods delivery where consumers can obtain or exchange faulty goods ordered through the online store. Online retailers wish to place the points of goods delivery in accordance with defined optimality principle. In this paper, a compromise solution is chosen as the optimality principle.

---

[1] malafeyevoa@mail.ru
[2] redinskich@yandex.ru
[3] rumyantsevnn@gmail.com

## 2. THE FORMALIZATION OF THE MODEL OF MANY-AGENT INTERACTION IN SOCIAL-TRADING NETWORK

The online retailer network $(N,k)$, where $N$ is a finite set of nodes, $k$ is the flow capacity function that corresponds to each edge $(x,y)$ of network $(N,k)$ with a non-negative number $k(x,y) \geq 0$ is considered. Let the producers of goods be located in some nodes of the network. The unit price for the online-retailers at various producers of goods consists of a set of numbers $L = (l_1,...,l_n)$. From the production points the goods are delivered to the points of goods delivery. In the network consumers of the goods are located in the set of $B = (b_1,...,b_j,...b_s)$. The value of the function of transport corruption costs of online-retailers $c_a(x,y)$ is determined for each edge of the network : $c_a(x,y) = c_a^{'}(x,y) + \gamma_a(x,y)$, where $c_a^{'}(x,y)$ is transportation costs, $\gamma_a(x,y)$ is corruption costs on the edge $(x, y)$ of the online- retailer. Also, for each network edge $(x, y)$ the value of the function of transport corruption costs $c_b(x,y)$ of the consumers (here under transport corruption costs of consumers we can understand selling to consumers expired, low-quality goods, goods not corresponding to its price, defective goods etc.) is determined: $c_b(x,y) = c_b^{'}(x,y) + \gamma_b(x,y)$, where $c_b^{'}(x,y)$ is transportation costs, $\gamma_b(x,y)$ is corruption costs on the edge $(x, y)$ for consumer. In the free nodes of the network there can be $h$ points of delivery $M = (m_1,...,m_i,...m_h)$, to each of which the goods are delivered. Cost of product $P$ at the point of goods delivery is equal to the sum of its value at the point of production where the good is purchased, transportation costs on delivery and markup of the online-retailers: $P_h = l_h + c_h(x,y) + T_h$, where $l_h$ is the cost of good for the online retailer in the production point, $c_h(x,y)$ is the amount of transport corruption costs on the path from production point to point of goods delivery, $T_h$ is the markup of online retailer. Online retailers wish to place the points of the goods delivery at the nodes in accordance with an optimality principle. In this paper a compromise solution is taken as the optimality principle. For finding the compromise solution it is necessary to know the payoff function for each player.

Players are online-retailers in this problem. Payoff function of each player is the value of the quantity of the product that each consumer purchases.

*The algorithm for finding a compromise solution.*

1) Let $Z$ be the set of possible solutions $z_s \in Z$, where $z_s = ((p_1^{i_1} m_1^{j_1}),...,(p_l^{i_l} m_n^{j_n}))$ are feasible nodes of the of goods delivery location for determined prices, and $s$ is an index, which renumbered all the valid solutions and $s = 1,...,\overline{s}$. $R_1(z_s),...,R_n(z_s)$ is the income in the $i$-th point of the goods delivery. Let us calculate the ideal vector $M = [M_1,...,M_n] = (\max R_1(z_s)),...,(\max R_n(z_s)) -$ which is the maximal value of income of the $j$-th points of goods delivery.

2) For each possible solution $z_s, s = \overline{1,s}$ let us calculate the matrix of residuals, which would be $A^* = (M_1 - R_1(z_s)),...,(M_n - R_n(z_s))$, где $s = \overline{1,s}$

3) In the resulting matrix of residuals for each solution $z_s, s = \overline{1,s}$ from the columns of matrix $A^*$ let us choose a maximal value $\gamma_s = \max_s\{(M_1 - R_1(z_s)),...,(M_n - R_n(z_s))\}$

4) Let us choose the minimum of these maximum solutions $\min_s = \gamma_s$, that will be the compromise solution.

## 3. EXAMPLES

*A. Compromise principle of optimality in network game host.*

The network $(N, k)$ (Fig. 1) containing 30 nodes and 51 edges is considered. Two values of functions of the transport corruption costs, one for online retailers and another one for the consumers are given on the edges of the network. The goods producers are located in the nodes $x_0$ and $x_{29}$. The unit value of the product for online-retailers in these vertices are $l_1 = 2$ and $l_2 = 3$ respectively. Four nodes $x_{22}, x_{15}, x_{12}, x_{16}$ are taken as feasible points of goods location delivery in the network. The points of delivery of the goods can not be located in one node. Markup of the online-retailer is 100% of the amount of the costs of the online-retailer for the purchase of goods and for the transport of goods along the way from point of production to point of goods delivery. The consumers are located in the nodes $x_8, x_{21}, x_{19}, x_{26}$.

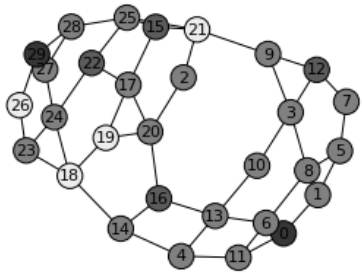

Fig 1. The network $(N, k)$

Transport corruption costs $c_a(x,y) = c'_a(x,y) + \gamma_a(x,y)$ for online retailers and $c_b(x,y) = c'_b(x,y) + \gamma_b(x,y)$ for consumers are given in table 1, where $(x, y)$ is the edge of the network, $c'_a(x,y)$ are transportation costs, $\gamma_a(x,y)$ are corruption costs on the edge $(x, y)$ for online-retailers, $c'_b(x,y)$ are transportation costs, $\gamma_b(x,y)$ are corruption costs on the edge $(x, y)$ for consumers:

TABLE I. TRANSPORT CORRUPTION COSTS

| (x,y) | $c'_a(x,y)$ | $\gamma_a(x,y)$ | $c'_b(x,y)$ | $\gamma_b(x,y)$ | $c_a(x,y)$ | $c_b(x,y)$ |
|---|---|---|---|---|---|---|
| (0,11) | 1 | 3 | 6 | 2 | 4 | 8 |
| (11,4) | 12 | 6 | 6 | 3 | 18 | 9 |
| (4,14) | 0 | 1 | 3 | 5 | 1 | 8 |
| (14,18) | 12 | 9 | 2 | 4 | 21 | 6 |
| (18,23) | 13 | 2 | 6 | 5 | 15 | 11 |
| (23,26) | 0 | 8 | 1 | 5 | 8 | 6 |
| (26,29) | 6 | 7 | 4 | 0 | 13 | 4 |
| (29,28) | 6 | 7 | 3 | 1 | 13 | 4 |
| (28,25) | 14 | 6 | 6 | 3 | 20 | 9 |
| (25,21) | 4 | 9 | 2 | 4 | 13 | 6 |
| (21,9) | 13 | 7 | 0 | 2 | 20 | 2 |
| (9,12) | 3 | 3 | 2 | 3 | 6 | 5 |
| (12,7) | 2 | 2 | 6 | 4 | 4 | 10 |
| (7,5) | 13 | 0 | 5 | 1 | 13 | 6 |
| (5,1) | 1 | 6 | 8 | 1 | 7 | 9 |
| (1,0) | 8 | 9 | 8 | 2 | 17 | 10 |
| (0,6) | 10 | 9 | 5 | 3 | 19 | 8 |
| (11,6) | 1 | 3 | 1 | 6 | 4 | 7 |
| (4,13) | 6 | 6 | 1 | 6 | 12 | 7 |
| (14,16) | 13 | 2 | 5 | 3 | 15 | 8 |
| (18,19) | 7 | 0 | 3 | 6 | 7 | 9 |
| (18,24) | 9 | 1 | 6 | 4 | 10 | 10 |
| (23,24) | 6 | 8 | 6 | 4 | 14 | 10 |
| (29,27) | 1 | 6 | 0 | 5 | 7 | 5 |
| (28,27) | 2 | 0 | 1 | 1 | 2 | 2 |
| (25,22) | 12 | 8 | 3 | 5 | 20 | 8 |
| (25,15) | 10 | 5 | 7 | 1 | 15 | 8 |
| (21,15) | 3 | 9 | 5 | 3 | 12 | 8 |
| (21,2) | 9 | 4 | 5 | 3 | 13 | 8 |
| (9,2) | 13 | 7 | 0 | 1 | 20 | 1 |
| (9,3) | 4 | 1 | 0 | 1 | 5 | 1 |
| (12,3) | 6 | 4 | 3 | 3 | 10 | 6 |
| (7,8) | 2 | 2 | 7 | 2 | 4 | 9 |
| (5,8) | 12 | 7 | 4 | 5 | 19 | 9 |
| (1,8) | 6 | 6 | 3 | 5 | 12 | 8 |
| (8,6) | 13 | 9 | 6 | 1 | 22 | 7 |
| (6,13) | 1 | 2 | 1 | 5 | 3 | 6 |
| (13,16) | 8 | 4 | 7 | 6 | 12 | 13 |
| (16,20) | 0 | 1 | 6 | 1 | 1 | 7 |
| (20,19) | 11 | 9 | 6 | 1 | 20 | 7 |
| (19,17) | 5 | 9 | 3 | 1 | 14 | 4 |
| (17,22) | 4 | 10 | 5 | 2 | 14 | 7 |
| (22,24) | 13 | 2 | 1 | 5 | 15 | 6 |
| (24,27) | 13 | 0 | 3 | 5 | 13 | 8 |
| (17,15) | 13 | 5 | 0 | 4 | 18 | 4 |
| (17,20) | 4 | 9 | 1 | 5 | 13 | 6 |
| (20,2) | 13 | 4 | 5 | 4 | 17 | 9 |
| (2,10) | 3 | 9 | 3 | 5 | 12 | 8 |
| (10,3) | 5 | 10 | 0 | 4 | 15 | 4 |
| (3,8) | 10 | 6 | 7 | 2 | 16 | 9 |
| (10,13) | 3 | 1 | 1 | 1 | 4 | 2 |

Let us calculate transport corruption costs for the four possible location points of goods delivery. Let us compute the weights of shortest paths in the graph from points of goods production to feasible points of goods delivery location:

TABLE II.

|       | $x_{22}$ | $x_{15}$ | $x_{12}$ | $x_{16}$ |
|-------|----------|----------|----------|----------|
| $x_0$ | 53 | 54 | 39 | 25 |
| $x_{29}$ | 38 | 47 | 71 | 61 |

Since the markup of the seller is 100%, the price of goods in all four feasible points of goods delivery location is:

TABLE III.

|       | $x_{22}$ | $x_{15}$ | $x_{12}$ | $x_{16}$ |
|-------|----------|----------|----------|----------|
| $x_0$ | 106 | 108 | 78 | 50 |
| $x_{29}$ | 76 | 94 | 142 | 122 |

Let us calculate the minimum cost of all paths from points of consumers location to four feasible points of goods delivery location:

TABLE IV.

|          | $x_{22}$ | $x_{15}$ | $x_{12}$ | $x_{16}$ |
|----------|----------|----------|----------|----------|
| $x_8$    | 26 | 20 | 15 | 26 |
| $x_{21}$ | 14 | 8  | 7  | 19 |
| $x_{19}$ | 11 | 8  | 22 | 14 |
| $x_{26}$ | 22 | 25 | 30 | 31 |

Let us calculate the income of the online retailers, depending on their choice of goods delivery location. There are only 6 situations:

TABLE V.

|   | (22,15) | (22,12) | (22,16) | (15,12) | (15,16) | (12,16) |
|---|---------|---------|---------|---------|---------|---------|
| 1 | 304 | 152 | 0   | 0   | 0   | 0   |
| 2 | 0   | 156 | 200 | 312 | 200 | 200 |

Let us $Z$ be the set of possible solutions $z_s \in Z$ where $z_s = ((p_1^{i_1} m_1^{j_1}),(p_1^{i_2} m_1^{j_2})), i_1, i_2 = 1,2; j_1, j_2 = 12,15,16,22$ are feasible points of goods delivery location for the given prices, and $s$ is an index, which renumbered all the valid solutions $s = 1,...,\bar{s}$. Let us calculate the ideal vector $M = [M_1,...,M_n] = [M_1, M_2] = (\max_s R_1(z_s), \max_s R_2(z_s))$, which is the maximal value of $j$-th online retailer: $M = (304, 312)$. For each possible solution $s = 1,...,\bar{s}$ let us compute the matrix of residuals, which would be $A^* = (M_1 - R_1(z_s)),(M_2 - R_2(z_s)) = (\alpha_s, \beta_s)$, where $s = 1,...,\bar{s}$. Let's denote $(\alpha_s, \beta_s)$ the first and second column correspondingly.

$$A^* = \begin{pmatrix} 0 & 312 \\ 152 & 156 \\ 304 & 112 \\ 304 & 0 \\ 304 & 112 \\ 304 & 112 \end{pmatrix}$$

In the matrix of residuals for each solution $z_s$, $s = 1,...,\bar{s}$, from the columns $(\alpha_s, \beta_s)$ let us choose maximal $\delta_s = \max_s\{\alpha_s, \beta_s\}$. Let us choose minimal from these maximal solutions $\min_s \delta_s = (152, 156)$, which will be the compromise solution.

According to the compromise solution, the points of goods delivery should be placed at points $x_{12}$ and $x_{22}$.

*B. Compromise-equilibrium principle of optimality allocation in the network game of two points of the goods delivery.*

The network $(N,k)$ (Fig. 2) containing 42 nodes and 70 edges is considered. On the edges of the network two values of transport corruption costs for online-retailers and consumers are given. In the nodes $x_0$ and $x_{12}$ the points of goods production are located. The cost of a unit of goods unit for online retailers in given nodes is $l_1 = 2$ and $l_2 = 3$ correspondingly.

Four nodes $x_5, x_{21}, x_{37}, x_{39}$ for feasible points of goods delivery location are selected in the network. The points of goods delivery cannot be located in one node. Markup of the online retailer is 100% of the amount of the costs of the online retailer for the purchase of goods and for the transportation of goods along the way from production point to point of goods delivery. The consumers are located in the nodes $x_8, x_9, x_{18}, x_{20}, x_{33}, x_{41}$.

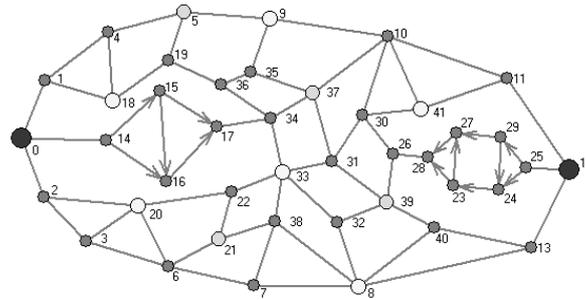

Fig 2. The network $(N,k)$

On the part of the network (0,34) and (12,26) transport corruption costs depend on the number of goods *n* which is necessary to bring the edge and is equal to *10n* for the edges (14,15),(16,17),(25,29) and (27,28), *20 + n* for the edges (15,17) and (14,16), *10 + n* for the edges (15,16),(25,24) and (23,28), *5 + n* for the edges (29,27) and (24,23), *n* for (29,24) and (27,23). Let us consider a section of the network (0,34). It consists of seven edges, where it is necessary to transport three units of goods from the node $x_0$ to the node $x_{34}$. Let us consider the 3 paths in which the goods can transfer, and we set up a system of equations to find out the solution corresponding to the Cournot-Nash equilibrium. Let us make a table with possible ways (TABLE VI).

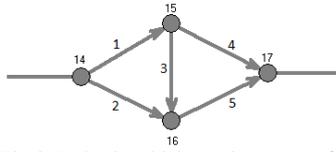

Fig 3. Paths, in which goods can transfer

TABLE VI.

| path's number | the numbers of edges included in the path |
|---|---|
| 1 | 1 – 4 |
| 2 | 1 - 3 – 5 |
| 3 | 2 – 5 |

Costs, required to cover each segment, depend on the number of goods in this segment. Let $x_1$, $x_2$, $x_3$ be the quantity of goods transported along the path with the corresponding index. Then, accordingly, the costs corresponding to each path can be written as:

1) $10(x_1 + x_2) + 20 + x_1$ – for the first path
2) $10(x_1 + x_2) + 10 + x_1 + 10(x_2 + x_3)$ – for the second path
3) $20 + x_3 + 10(x_2 + x_3)$ – for the third path

After simplification we get the following expressions:

1) $11x_1 + 10x_2 + 20$
2) $10x_1 + 21x_2 + 10x_3 + 10$
3) $10x_2 + 11x_3 + 20$

For the Cournot-Nash equilibrium all of these values must be equal, this leads to the system of equations:

$11x_1 + 10x_2 + 20 = 10x_1 + 21x_2 + 10x_3 + 10 = 10x_2 + 11x_3 + 20$ We know, that $x_1 + x_2 + x_3 = 3$, then we can solve the following system of equations:

$$\begin{cases} x_1 - 11x_2 - 10x_3 = -10 \\ 11x_1 - 11x_3 = 0 \\ x_1 + x_2 + x_3 = 3 \end{cases}$$

The solution of this system is:

$$\begin{cases} x_1 = \frac{23}{13} \\ x_2 = -\frac{7}{13} \\ x_3 = \frac{23}{13} \end{cases}$$

By substituting the obtained solution to the original system we calculate the total cost of the path (0,34), which will be equal $34\frac{1}{13}$. Next, let us consider the section of the network (12,26) consisting of 10 edges, where it is necessary to transfer three items from the node $x_{12}$ to the node $x_{26}$. Let us consider 4 paths, in which the goods can transfer, and we set up a system of equations to find out the solution corresponding to the Cournot-Nash equilibrium.

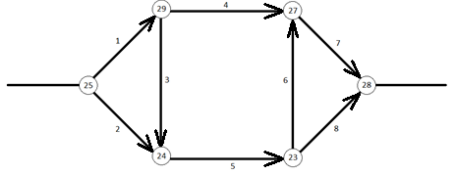

Fig 4. Paths, in which goods can transfer

Let us make a table that will describe possible paths:

TABLE VII.

| path's number | the numbers of edges included in the path |
|---|---|
| 1 | 1 - 4 – 7 |
| 2 | 1 - 3 - 5 - 6 – 7 |
| 3 | 2 - 5 - 6 - 7 |
| 4 | 2 - 5 – 8 |

Costs, required to cover each segment, depend on the number of goods in this segment. Let $x_1, x_2, x_3, x_4$ be the quantity of goods transported along the path with the corresponding index. Then, accordingly, the costs corresponding to each path can be written as:

1) $10(x_1 + x_2) + 5 + x_1 + 10(x_1 + x_2 + x_3)$ – for the first path
2) $10(x_1 + x_2) + 5 + x_1 + (x_1 + x_2 + x_3) + x_2 + x_3 +$
   $+ 10(x_1 + x_2 + x_3)$ – for the second path
3) $10 + (x_3 + x_4) + 5 + (x_2 + x_3 + x_4) + x_2 + x_3 +$
   $+ 10(x_1 + x_2 + x_3)$ – for the third path
4) $10 + (x_3 + x_4) + 5 + (x_2 + x_3 + x_4) + 10 + x_4$ – for the fourth path

After simplification we get the following expressions:
1) $21x_1 + 20x_2 + 10x_3 + 5$
2) $20x_1 + 23x_2 + 12x_3 + x_4 + 5$
3) $10x_1 + 12x_2 + 13x_3 + 2 + x_4 + 15$
4) $x_2 + 2x_3 + 3x_4 + 25$

For the Cournot-Nash equilibrium all of these values must be equal. That leads to the system of equations:

$21x_1 + 20x_2 + 10x_3 + 5 = 20x_1 + 23x_2 + 12x_3 + x_4 + 5 =$
$= 10x_1 + 12x_2 + 13 + 2x_4 + 15 = x_2 + 2x_3 + 3x_4 + 25$

We know that $x_1 + x_2 + x_3 + x_4 = 3$, then we can solve the following system of equations:

$$\begin{cases} x_1 - 3x_2 - 2x_3 - x_4 = 0 \\ 11x_1 + 8x_2 - 3x_3 - 2x_4 = 10 \\ 21x_1 + 19x_2 + 8x_3 - 3x_4 = 20 \\ x_1 + x_2 + x_3 + x_4 = 3 \end{cases}$$

The solution of this system is:

$$\begin{cases} x_1 = \dfrac{31}{23} \\ x_2 = -\dfrac{7}{46} \\ x_3 = 0 \\ x_4 = \dfrac{83}{46} \end{cases}$$

By substituting the obtained solution to the original system we calculate the total cost of the path (12,26), which will be equal $30\dfrac{6}{23}$. Transport corruption costs $c_a(x,y) = c_a'(x,y) + \gamma_a(x,y)$ for the online retailers and $c_b(x,y) = c_b'(x,y) + \gamma_b(x,y)$ for the consumers are shown in table 3, where $(x, y)$ is the edge of the network, $c_a'(x, y)$ are transportation costs, $\gamma_a(x, y)$ are corruption costs on the edge $(x, y)$ for consumers, $c_b'(x, y)$ are transportation costs, $\gamma_b(x, y)$ are corruption costs on the edge $(x, y)$ for consumers:

TABLE VIII. Transport corruption costs

| (x,y) | $c_a'(x,y)$ | $\gamma_a(x,y)$ | $c_b'(x,y)$ | $\gamma_b(x,y)$ | $c_a(x,y)$ | $c_b(x,y)$ |
|---|---|---|---|---|---|---|
| (0,1) | 15 | 15 | 15 | 15 | 30 | 30 |
| (1,4) | 4 | 9 | 4 | 5 | 13 | 9 |
| (4,5) | 8 | 0 | 5 | 0 | 8 | 5 |
| (5,9) | 14 | 2 | 3 | 6 | 16 | 9 |
| (9,10) | 1 | 6 | 1 | 2 | 7 | 3 |
| (10,11) | 5 | 7 | 3 | 3 | 12 | 6 |
| (11,12) | 15 | 15 | 15 | 15 | 30 | 30 |
| (12,13) | 15 | 15 | 15 | 15 | 30 | 30 |
| (13,8) | 3 | 1 | 3 | 3 | 4 | 6 |
| (8,7) | 4 | 5 | 2 | 5 | 9 | 7 |
| (7,6) | 12 | 4 | 2 | 2 | 16 | 4 |
| (6,3) | 12 | 5 | 1 | 6 | 17 | 7 |
| (3,2) | 6 | 1 | 4 | 1 | 7 | 5 |
| (2,0) | 15 | 15 | 15 | 15 | 30 | 30 |
| (1,18) | 3 | 3 | 2 | 6 | 6 | 8 |
| (4,18) | 8 | 2 | 7 | 6 | 10 | 13 |
| (5,19) | 7 | 3 | 3 | 2 | 10 | 5 |
| (9,35) | 6 | 1 | 3 | 3 | 7 | 6 |
| (10,37) | 7 | 10 | 3 | 3 | 17 | 6 |
| (10,30) | 9 | 9 | 7 | 2 | 18 | 9 |
| (10,41) | 12 | 6 | 8 | 5 | 18 | 13 |
| (11,41) | 2 | 10 | 5 | 5 | 12 | 10 |

| | | | | | |
|---|---|---|---|---|---|
| (13,40) | 6 | 10 | 2 | 0 | 16 | 2 |
| (8,40) | 13 | 6 | 4 | 6 | 19 | 10 |
| (8,32) | 14 | 5 | 4 | 3 | 19 | 7 |
| (8,38) | 8 | 2 | 6 | 5 | 10 | 11 |
| (7,38) | 13 | 4 | 1 | 2 | 17 | 3 |
| (6,21) | 5 | 1 | 2 | 2 | 6 | 4 |
| (6,20) | 5 | 3 | 6 | 3 | 8 | 9 |
| (3,20) | 6 | 2 | 4 | 1 | 8 | 5 |
| (2,20) | 1 | 4 | 0 | 5 | 5 | 5 |
| (20,22) | 4 | 7 | 1 | 1 | 11 | 2 |
| (22,21) | 13 | 0 | 2 | 4 | 13 | 6 |
| (21,38) | 1 | 8 | 2 | 4 | 9 | 6 |
| (38,33) | 5 | 3 | 4 | 2 | 8 | 6 |
| (22,33) | 4 | 6 | 1 | 6 | 10 | 7 |
| (33,31) | 13 | 4 | 1 | 5 | 17 | 6 |
| (18,19) | 13 | 3 | 8 | 3 | 16 | 11 |
| (19,36) | 13 | 3 | 3 | 0 | 16 | 3 |
| (36,35) | 5 | 3 | 5 | 4 | 8 | 9 |
| (35,37) | 10 | 4 | 3 | 3 | 14 | 6 |
| (37,34) | 2 | 7 | 5 | 4 | 9 | 9 |
| (34,36) | 5 | 6 | 0 | 1 | 11 | 1 |
| (34,33) | 12 | 9 | 1 | 6 | 21 | 7 |
| (31,37) | 11 | 6 | 1 | 0 | 17 | 1 |
| (33,32) | 3 | 6 | 1 | 4 | 9 | 5 |
| (32,39) | 2 | 6 | 3 | 2 | 8 | 5 |
| (39,31) | 2 | 6 | 2 | 5 | 8 | 7 |
| (39,40) | 7 | 7 | 4 | 4 | 14 | 8 |
| (39,26) | 3 | 0 | 5 | 3 | 3 | 8 |
| (26,30) | 4 | 9 | 2 | 1 | 13 | 3 |
| (41,30) | 5 | 2 | 3 | 1 | 7 | 4 |
| (30,31) | 13 | 0 | 2 | 0 | 13 | 2 |

Let us calculate the transport corruption costs in the four feasible points of goods delivery location. We use the Floyd algorithm and calculate the weights of shortest paths in the graph from points of production to feasible points of goods delivery location:

TABLE IX.

| | $x_5$ | $x_{21}$ | $x_{37}$ | $x_{39}$ |
|---|---|---|---|---|
| $x_0$ | 71 | 72 | 43 | 68 |
| $x_{12}$ | 84 | 67 | 58 | 33 |

Since the markup for the online-retailer is 100%, the price of goods in all four feasible points of goods delivery location is:

TABLE X.

|          | $x_5$ | $x_{21}$ | $x_{37}$ | $x_{39}$ |
|----------|-------|----------|----------|----------|
| $x_0$    | 146   | 148      | 90       | 140      |
| $x_{12}$ | 174   | 140      | 122      | 72       |

Next, let us calculate the matrix of minimum of weights of all paths from points of location of all consumers in all feasible points of goods delivery location:

TABLE XI.

|          | $x_5$ | $x_{21}$ | $x_{37}$ | $x_{39}$ |
|----------|-------|----------|----------|----------|
| $x_8$    | 28    | 15       | 19       | 12       |
| $x_9$    | 9     | 28       | 9        | 17       |
| $x_{18}$ | 16    | 34       | 24       | 32       |
| $x_{20}$ | 25    | 8        | 16       | 19       |
| $x_{33}$ | 16    | 12       | 7        | 10       |
| $x_{41}$ | 25    | 24       | 7        | 13       |

Let us calculate the matrix of income of the online retailers, depending on their choice of feasible points of goods delivery location. We have a total of 6 situations:

TABLE XII.

|   | (5,21) | (5,37) | (5,39) | (21,37) | (21,39) | (37,39) |
|---|--------|--------|--------|---------|---------|---------|
| 1 | 292    | 0      | 0      | 0       | 0       | 0       |
| 2 | 560    | 540    | 432    | 540     | 432     | 432     |

Let us $Z$ be a set of possible locations $z_s \in Z$ where $z_s = ((p_1^{i_1} m_1^{j_1}),(p_1^{i_2} m_1^{j_2})), i_1, i_2 = 1,2; j_1, j_2 = 5,21,37,39$ are possible nodes of goods delivery locations by given price, and $s$ is an index, that renumbered all possible solutions $s = 1,...,\bar{s}$. Let us compute the ideal vector $M = [M_1,...,M_n] = [M_1, M_2] = (\max_s R_1(z_s), \max_s R_2(z_s))$ which is the maximal value of income of $j$-th online retailer: $M = (292, 560)$. For each possible solution $s = 1,...,\bar{s}$ let us calculate the matrix of residuals, that will be $A^* = (M_1 - R_1(z_s)),(M_2 - R_2(z_s)) = (\alpha_s, \beta_s)$, where $s = 1,...,\bar{s}$. Let us denote $(\alpha_s, \beta_s)$ the first and the second column correspondingly.

$$A^* = \begin{pmatrix} 0 & 0 \\ 292 & 20 \\ 292 & 128 \\ 292 & 20 \\ 292 & 128 \\ 292 & 128 \end{pmatrix}$$

In the resulting matrix of residuals for each solution $z_s$, $s = 1,...,\bar{s}$, from the columns $(\alpha_s, \beta_s)$ we choose maximal $\delta_s = \max_s \{\alpha_s, \beta_s\}$. Let us choose the minimal from these maximal solutions $\min_s \delta_s = (0,0)$ that will be compromise solution.

According to the compromise solution, the points of goods delivery should be placed at points $x_5$ and $x_{21}$.

## 4. CONCLUSION

In this paper the model of online-retailers network is formalized. The algorithm for finding the optimal choosing of goods delivery location in accordance with the compromise solution as the principle of optimality is proposed. Two numerical examples are solved.

## 5. ACKNOWLEDGMENT

The work was supported by RFBR, project No. 14-06-00326****.